# Some Results and Characterizations for Mannheim Offsets of the Ruled Surfaces

## Mehmet ÖNDER[*], H. Hüseyin UĞURLU[**]

[*]Celal Bayar University, Faculty of Science and Arts, Department of Mathematics, Muradiye Manisa. E-mail: mehmet.onder@bayar.edu.tr

[**]Gazi University, Faculty of Education, Department of Secondary Education Science and Mathematics Teaching, Mathematics Teaching Programme, Ankara, Turkey
E-mail: hugurlu@gazi.edu.tr

**Abstract**

In this study, we give the dual characterizations of Mannheim offsets of the ruled surface in terms of their integral invariants and the new characterization of the Mannheim offsets of developable surface. Furthermore, we obtain the relationships between the area of projections of spherical images for Mannheim offsets of ruled surfaces and their integral invariants.

## 1. Introduction

Ruled surfaces are the surfaces which are generated by moving a straight line continuously in the space and are one of the most important topics of differential geometry. A ruled surface can always be described (at least locally) as the set of points swept by a moving straight line. These surfaces are used in many areas of sciences such as Computer Aided Geometric Design (CAGD), mathematical physics, moving geometry, kinematics for modeling the problems and model-based manufacturing of mechanical products. Especially, the offsets of ruled surface have an important role in (CAGD). Some studies dealing with offsets of the surfaces have been given in[2,8,10,12,13,14]. In [14], Ravani and Ku have defined and given a generalization of the theory of Bertrand curves for Bertrand trajectory ruled surfaces on the line geometry. Küçük and Gürsoy have studied the integral invariants of Bertrand trajectory ruled surfaces in dual space and given the relations between the invariants[8].

Furthermore, similar to the Bertrand curves, in [9] a new definition of special curve pair has been given by Liu and Wang: Let $C$ and $C^*$ be two space curves. $C$ is said to be a Mannheim partner curve of $C^*$ if there exists a one to one correspondence between their points such that the binormal vector of $C$ is the principal normal vector of $C^*$. Orbay and et.al. have given a generalization of the theory of Mannheim curves for ruled surfaces and called Mannheim offsets[10].

In this paper, we examine the Mannheim offsets of trajectory ruled surfaces in view of their integral invariants. Using the dual representations of the ruled surfaces, we give the results obtained in [10] and some new results in short forms. Moreover, we show that if the Mannheim offsets of trajectory ruled surfaces are developable then their striction lines are Mannheim partner curves. Furthermore, we give some characterizations of Mannheim offsets of trajectory ruled surfaces in terms of integral invariants (such as the angle of pitch and the pitch) of closed trajectory ruled surfaces. Finally, we obtain the relationship between the area of projections of spherical images of Mannheim offsets of trajectory ruled surfaces and their integral invariants.

## 2. Differential Geometry of Ruled Surfaces in $E^3$

Let $I$ be an open interval in the real line $IR$. Let $k = k(s)$ be a curve in $E^3$ defined on $I$ and $q = q(s)$ be a unit direction vector of an oriented line in $IR^3$. Then we have the following parametrization for a ruled surface

$$\varphi_q(s,v) = \vec{k}(s) + v\vec{q}(s). \tag{1}$$

The parametric $s$-curve of this surface is a straight line of the surface which is called ruling. For $v=0$, the parametric $v$-curve of this surface is $\vec{k} = \vec{k}(s)$ which is called base curve or generating curve of the surface. In particular, if $\vec{q}$ is constant, the ruled surface is said to be cylindrical, and non-cylindrical otherwise[7].

The striction point on a ruled surface is the foot of the common normal between two consecutive rulings. The set of the striction points constitute a curve $\vec{c} = \vec{c}(s)$ lying on the ruled surface and is called striction curve. The parametrization of the striction curve $\vec{c} = \vec{c}(s)$ on a ruled surface is given by

$$\vec{c}(s) = \vec{k}(s) - \frac{\langle d\vec{q}, d\vec{k} \rangle}{\langle d\vec{q}, d\vec{q} \rangle} \vec{q}. \tag{2}$$

So that, the base curve of the ruled surface is its striction curve if and only if $\langle d\vec{q}, d\vec{k} \rangle = 0$ [7].

The distribution parameter (or drall) of the ruled surface in (1) is given as

$$\delta_q = \frac{\langle d\vec{k}, \vec{q} \times d\vec{q} \rangle}{\langle d\vec{q}, d\vec{q} \rangle} \tag{3}$$

If $\delta_q = 0$, then the normal vectors of the ruled surface are collinear at all points of the same ruling and at the nonsingular points of the ruled surface the tangent planes are identical. We then say that the tangent plane contacts the surface along a ruling. Such a ruling is called a *torsal* ruling. If $\delta_q \neq 0$, then the tangent planes are distinct at all points of the same ruling which is called *nontorsal*.

A ruled surface whose all rulings are torsal is called a *developable ruled surface*. The remaining ruled surfaces are called *skew ruled surfaces*. Thus, from (3) a ruled surface is developable if and only if at all its points the distribution parameter $\delta_q = 0$ [7,14].

Let $\left\{ \vec{q}, \vec{h} = \dfrac{d\vec{q}/ds}{\|d\vec{q}/ds\|}, \vec{a} = \vec{q} \times \vec{h} \right\}$ be a moving othonormal trihedron making a spatial motion along a closed space curve $\vec{k}(s)$, $s \in \mathbb{R}$, in $E^3$. In this motion, an oriented line fixed in the moving system generates a closed ruled surface called closed trajectory ruled surface (CTRS) in $E^3$ [8]. A parametric equation of a closed trajectory ruled surface generated by $\vec{q}$-axis is

$$\varphi_q(s,v) = \vec{k}(s) + v\vec{q}(s), \quad \varphi(s+2\pi,v) = \varphi(s,v), \quad s,v \in \mathbb{R}. \tag{4}$$

Consider the moving orthonormal system $\{\vec{q}, \vec{h}, \vec{a}\}$. Then, the axes of the trihedron intersect at the striction point of $\vec{q}$-generator of $\varphi_q$-CTRS. The structral equations of this motion are

$$\begin{cases} d\vec{q} = k_1 \vec{h} \\ d\vec{h} = -k_1 \vec{q} + k_2 \vec{a} \\ d\vec{a} = -k_2 \vec{h} \end{cases} \tag{5}$$

and

$$\frac{db}{ds} = \cos\sigma \vec{q} + \sin\sigma \vec{a} \tag{6}$$

where $b = b(s)$ is the striction line of $\varphi_q$-CTRS and the differential forms $k_1$, $k_2$ and $\sigma$ are the natural curvature, the natural torsion and the striction of $\varphi_q$-CTRS, respectively[7,8]. Here, the striction is restricted as $-\pi/2 < \sigma < \pi/2$ for the orientation on $\varphi_q$-CTRS and $s$ is the length of the striction line.

The pole vector and the Steiner vector of the motion are given by

$$\vec{p} = \frac{\vec{\psi}}{\|\vec{\psi}\|}, \quad \vec{d} = \oint \vec{\psi} \tag{7}$$

respectively, where $\vec{\psi} = k_2 \vec{q} + k_1 \vec{a}$ is the instantaneous Pfaffian vector of the motion.

The pitch of $\varphi_q$-CTRS is defined by

$$\ell_q = \oint d\mu = -\oint \langle d\vec{k}, \vec{q} \rangle. \tag{8}$$

The angle of pitch of $\varphi_q$-CTRS is given one of the followings

$$\lambda_q = \oint d\theta = -\oint \langle d\vec{h}, \vec{a} \rangle = -\langle \vec{q}, \vec{d} \rangle = 2\pi - a_q = \oint g_q \tag{9}$$

where $a_q$ and $g_q$ are the measure of the spherical surface area bounded by the spherical image of $\varphi_q$-CTRS and the geodesic curvature of this image, respectively. The pitch and the angle of pitch are well-known real integral invariants of closed trajectory ruled surface[3-6].

The area vector of a $x$-closed space curve in $E^3$ is given by

$$v_x = \oint x \times dx \tag{10}$$

and the area of projection of a $x$-closed space curve in direction of the generator of a $y$-CTRS is

$$2 f_{x,y} = \langle v_x, y \rangle. \tag{11}$$

(See [5]).

## 3. Dual Numbers and Dual Vectors

Dual numbers had been introduced by W. K. Clifford (1845-1879). A dual number has the form $\bar{a} = (a, a^*) = a + \varepsilon a^*$ where $a$ and $a^*$ are real numbers and $\varepsilon = (0,1)$ is dual unit with $\varepsilon^2 = 0$. The product of dual numbers $\bar{a} = (a, a^*) = a + \varepsilon a^*$ and $\bar{b} = (b, b^*) = b + \varepsilon b^*$ is given by

$$\bar{a}\bar{b} = (a, a^*)(b, b^*) = (ab, ab^* + a^*b) = ab + \varepsilon(ab^* + a^*b). \tag{12}$$

We denote the set of dual numbers by $D$:

$$D = \{\bar{a} = a + \varepsilon a^* : a, a^* \in IR, \varepsilon^2 = 0\}. \tag{13}$$

Clifford showed that dual numbers form algebra, but not a field. The pure dual numbers $\varepsilon a^*$ are zero divisors, $(\varepsilon a^*)(\varepsilon b^*) = 0$. However, the other laws of the algebra of dual numbers are the same as the laws of algebra of complex numbers. This means that dual numbers form a ring over the real number field.

Now let $f$ be a differentiable function with dual variable $\bar{x} = x + \varepsilon x^*$. Then the Maclaurine series generated by $f$ is given by

$$f(\bar{x}) = f(x + \varepsilon x^*) = f(x) + \varepsilon x^* f'(x), \tag{14}$$

where $f'(x)$ is derivative of $f(x)$.

Let $D^3$ be the set of all triples of dual numbers, i.e.,

$$D^3 = \{\tilde{a} = (\bar{a}_1, \bar{a}_2, \bar{a}_3) : \bar{a}_i \in D, i = 1,2,3\}, \tag{15}$$

Then the set $D^3$ is called dual space. The elements of $D^3$ are called dual vectors. A dual vector $\tilde{a}$ may be expressed in the form $\tilde{a} = \vec{a} + \varepsilon \vec{a}^* = (\vec{a}, \vec{a}^*)$, where $\vec{a}$ and $\vec{a}^*$ are the vectors of $IR^3$. Then for any vectors $\tilde{a} = \vec{a} + \varepsilon \vec{a}^*$ and $\tilde{b} = \vec{b} + \varepsilon \vec{b}^*$ in $D^3$, the scalar product and the vector product are defined by

$$\langle \tilde{a}, \tilde{b} \rangle = \langle \vec{a}, \vec{b} \rangle + \varepsilon \left( \langle \vec{a}, \vec{b}^* \rangle + \langle \vec{a}^*, \vec{b} \rangle \right), \tag{16}$$

and

$$\tilde{a} \times \tilde{b} = \vec{a} \times \vec{b} + \varepsilon \left( \vec{a} \times \vec{b}^* + \vec{a}^* \times \vec{b} \right), \tag{17}$$

respectively, where $\langle \vec{a}, \vec{b} \rangle$ and $\vec{a} \times \vec{b}$ are the inner product and the vector product of the vectors where $\vec{a}$ and $\vec{a}^*$ in $IR^3$.

The norm of a dual vector $\tilde{a}$ is given by

$$\|\tilde{a}\| = \sqrt{\langle \tilde{a}, \tilde{a} \rangle} = \|\vec{a}\| + \varepsilon \frac{\langle \vec{a}, \vec{b}^* \rangle}{\|\vec{a}^*\|}. \tag{18}$$

A dual vector $\tilde{a}$ with norm 1 is called dual unit vector. The set of dual unit vectors is

$$\tilde{S}^2 = \left\{ \tilde{a} = (a_1, a_2, a_3) \in D^3 : \langle \tilde{a}, \tilde{a} \rangle = 1 \right\}, \tag{19}$$

which is called dual unit sphere.

E. Study used dual numbers and dual vectors in his research on the geometry of lines and kinematics. He devoted special attention to the representation of directed lines by dual unit vectors and defined the mapping that is known by his name: There exists one-to-one correspondence between the vectors of dual unit sphere $\tilde{S}^2$ and the directed lines of space of lines $IR^3$. By the aid of this correspondence, the properties of the spatial motion of a line can be derived. Hence, the geometry of ruled surface is represented by the geometry of dual curves on the dual unit sphere in $D^3$ (Fig. 1).

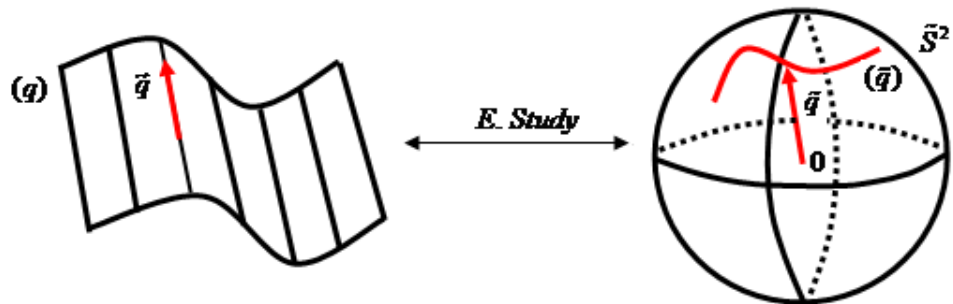

**Fig. 1.** E. Study Mapping

The angle $\bar{\theta} = \theta + \varepsilon \theta^*$ between two dual unit vectors $\tilde{a}, \tilde{b}$ is called *dual angle* and defined by

$$\langle \tilde{a}, \tilde{b} \rangle = \cos \bar{\theta} = \cos \theta - \varepsilon \theta^* \sin \theta. \tag{20}$$

By considering The E. Study Mapping, the geometric interpretation of dual angle is that $\theta$ is the real angle between lines $L_1$, $L_2$ corresponding to the dual unit vectors $\tilde{a}$, $\tilde{b}$, respectively, and $\theta^*$ is the shortest distance between those lines[1,6].

Let $K$ be a moving dual unit sphere generated by a dual orthonormal system

$$\left\{ \tilde{q},\ \tilde{h} = \frac{d\tilde{q}}{\|d\tilde{q}\|},\ \tilde{a} = \tilde{q} \times \tilde{h} \right\},\ \tilde{q} = \vec{q} + \varepsilon \vec{q}^*,\ \tilde{h} = \vec{h} + \varepsilon \vec{h}^*,\ \tilde{a} = \vec{a} + \varepsilon \vec{a}^*, \qquad (21)$$

and $K'$ be a fixed dual unit sphere with the same center. Then, the derivative equations of the dual spherical closed motion of $K$ with respect to $K'$ are

$$\begin{cases} d\tilde{q} = \bar{k}_1 \tilde{h} \\ d\tilde{h} = -\bar{k}_1 \tilde{q} + \bar{k}_2 \tilde{a} \\ d\tilde{a} = -\bar{k}_2 \tilde{h} \end{cases} \qquad (22)$$

where $\bar{k}_1(s) = k_1(s) + \varepsilon k_1^*(s)$, $\bar{k}_2(s) = k_2(s) + \varepsilon k_2^*(s)$, $(s \in \mathbb{R})$ are dual curvature and dual torsion, respectively. From the E. Study mapping, during the spherical motion of $K$ with respect to $K'$, the dual unit vector $\tilde{q}$ draws a dual curve on dual unit sphere $K'$ and this curve represents a ruled surface with ruling $\vec{q}$ in line space $\mathbb{R}^3$.

Dual vector $\tilde{\psi} = \vec{\psi} + \varepsilon \vec{\psi}^* = \bar{k}_2 \tilde{q} + \bar{k}_1 \tilde{a}$ is called the instantaneous Pfaffian vector of the motion and the vector $\tilde{P}$ given by $\tilde{\psi} = \|\tilde{\psi}\| \tilde{P}$ is called the dual pole vector of the motion. Then the vector

$$\tilde{d} = \oint \tilde{\psi} \qquad (23)$$

is the dual Steiner vector of the closed motion.

By considering the E. Study mapping, the dual equations (22) correspond to the real equations (5) and (6) of a closed spatial motion in $E^3$. So, the differentiable dual closed curve $\tilde{q} = \tilde{q}(s)$ is corresponds to a closed trajectory ruled surface in line space and denoted by $\varphi_q$-CTRS.

A dual integral invariant of a $\varphi_q$-CTRS can be given in terms of real integral invariants as follows and is called the dual angle of pitch of a $\varphi_q$-CTRS

$$\bar{\lambda}_q = -\oint \langle d\tilde{h}, \tilde{a} \rangle = -\langle \tilde{q}, \tilde{d} \rangle = 2\pi - \bar{a}_q = \oint \bar{g}_q = \lambda_q - \varepsilon \ell_q \qquad (24)$$

where $\tilde{d} = \vec{d} + \varepsilon \vec{d}^*$, $\bar{a}_q = a_q + \varepsilon a_q^*$ and $\bar{g}_q = g_q + \varepsilon g_q^*$ are the dual Steiner vector of the motion, the measure of dual spherical surface area and the dual geodesic curvature of spherical image of $\varphi_q$-CTRS, respectively[1,3,6].

**4. Mannheim Offsets of Trajectory Ruled Surfaces**

Let $\varphi_q$ and $\varphi_{q_1}$ be two trajectory ruled surfaces generated by dual vectors $\tilde{q}$ and $\tilde{q}_1$ of the dual orthonormal frames $\{\tilde{q}(s),\ \tilde{h}(s),\ \tilde{a}(s)\}$ and $\{\tilde{q}_1(s_1),\ \tilde{h}_1(s_1),\ \tilde{a}_1(s_1)\}$, respectively. Then $\varphi_q$ and $\varphi_{q_1}$ are called Mannheim offsets of trajectory ruled surface, if

$$\tilde{a}(s) = \tilde{h}_1(s_1) \qquad (25)$$

where $s$ and $s_1$ are the arc-length of striction lines of $\varphi_q$ and $\varphi_{q_1}$, respectively. By this definition, the relation between the trihedrons

$$\left\{ \tilde{q},\ \tilde{h} = \frac{d\tilde{q}}{\|d\tilde{q}\|},\ \tilde{a} = \tilde{q} \times \tilde{h} \right\} \qquad (26)$$

and

$$\left\{ \tilde{q}_1,\ \tilde{h}_1 = \frac{d\tilde{q}_1}{\|d\tilde{q}_1\|},\ \tilde{a}_1 = \tilde{q}_1 \times \tilde{h}_1 \right\} \qquad (27)$$

of trajectory ruled surfaces $\varphi_q$ and $\varphi_{q_1}$ can be given as follows

$$\begin{pmatrix} \tilde{q}_1 \\ \tilde{h}_1 \\ \tilde{a}_1 \end{pmatrix} = \begin{pmatrix} \cos\bar{\theta} & \sin\bar{\theta} & 0 \\ 0 & 0 & 1 \\ \sin\bar{\theta} & -\cos\bar{\theta} & 0 \end{pmatrix} \begin{pmatrix} \tilde{q} \\ \tilde{h} \\ \tilde{a} \end{pmatrix} \qquad (28)$$

where $\bar{\theta} = \theta + \varepsilon\theta^*$, $(0 \leq \theta \leq \pi,\ \theta^* \in \mathbb{R})$ is the dual angle between the generators $\tilde{q}$ and $\tilde{q}_1$ of Mannheim trajectory ruled surface $\varphi_q$ and $\varphi_{q_1}$. The angle $\theta$ is called the offset angle and $\theta^*$ is called the offset distance(Fig. 2).

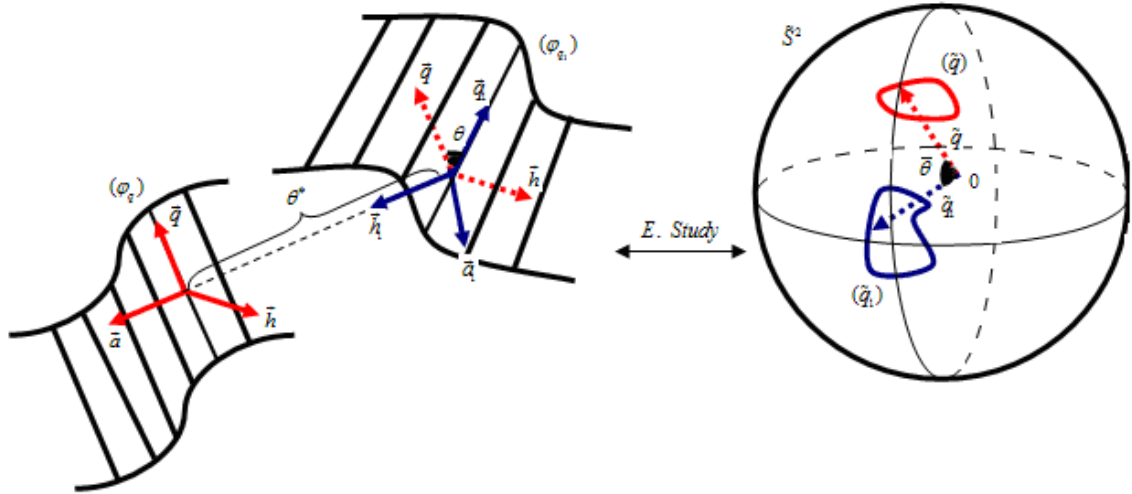

**Fig. 2.** E. Study Mapping of the Mannheim offset surfaces

Then, $\bar{\theta} = \theta + \varepsilon\theta^*$ is called dual offset angle of the Mannheim trajectory ruled surface $\varphi_q$ and $\varphi_{q_1}$. If $\theta = 0$ and $\theta = \pi/2$ then the Mannheim offsets are said to be oriented offsets and right offsets, respectively. Thus, we can give the followings.

***Theorem 4.1.*** *Let $\varphi_q$ and $\varphi_{q_1}$ be the Mannheim trajectory ruled surface. The offset angle and the offset distance are given by*

$$\theta = -\int k_1 ds, \quad \theta^* = -\int k_1^* ds \qquad (29)$$

*respectively.*
**Proof:** From (28) we have

$$\tilde{q}_1 = \cos\bar{\theta}\,\tilde{q} + \sin\bar{\theta}\,\tilde{h} \qquad (30)$$

Differentiating (30) and by using (22) and (28) we may write

$$\frac{d\tilde{q}_1}{ds} = -\left(\frac{d\bar{\theta}}{ds} + \bar{k}_1\right)\tilde{a}_1 + \bar{k}_2 \sin\bar{\theta}\,\tilde{h}_1 \qquad (31)$$

Since $\dfrac{d\tilde{q}_1}{ds}$ is orthogonal to $\tilde{a}_1$, from (31) we get

$$\bar{\theta} = -\int \bar{k}_1 ds \qquad (32)$$

Separating (32) into real and dual parts we have

$$\theta = -\int k_1 ds, \quad \theta^* = -\int k_1^* ds\,.$$

**Theorem 4.2.** *The closed trajectory ruled surfaces $\varphi_q$ and $\varphi_{q_1}$ form a Mannheim offsets if and only if the following relationship holds*

$$\overline{\wedge}_{q_1} = \overline{\wedge}_q \cos\overline{\theta} + \overline{\wedge}_h \sin\overline{\theta}, \quad \overline{\theta} = -\int \overline{k}_1 ds \tag{33}$$

**Proof:** Let the closed trajectory ruled surfaces $\varphi_q$ and $\varphi_{q_1}$ form a Mannheim offsets. Then, by direct calculation, from (24) and (28), the dual angle of pitch of $\varphi_{q_1}$-CTRS is given by

$$\overline{\wedge}_{q_1} = \overline{\wedge}_q \cos\overline{\theta} + \overline{\wedge}_h \sin\overline{\theta}$$

Conversely, if (33) holds, it is easily seen that $\varphi_q$ and $\varphi_{q_1}$-CTRS form a Mannheim offsets.

Equality (33) is a dual characterization of Mannheim offsets of CTRS in terms of their dual integral invariants. Separating (33) into real and dual parts, we obtain

$$\begin{cases} \lambda_{q_1} = \lambda_q \cos\theta + \lambda_h \sin\theta \\ \ell_{q_1} = (\ell_q - \theta^*\lambda_h)\cos\theta + (\ell_h + \theta^*\lambda_q)\sin\theta \end{cases} \tag{34}$$

Then, we may give the following results

**Result 4.1.** If $\varphi_q$ and $\varphi_{q_1}$ are the oriented closed Mannheim trajectory ruled surfaces i.e., $\theta = 0$, then the relationships between the real integral invariants of $\varphi_q$ and $\varphi_{q_1}$-CTRS are given as follows.

$$\lambda_{q_1} = \lambda_q, \quad \ell_{q_1} = \ell_q - \theta^*\lambda_h. \tag{35}$$

Furthermore, the measure of spherical surface areas bounded by the spherical images of $\varphi_q$ and $\varphi_{q_1}$-CTRS Mannheim offsets are the same, i.e.,

$$a_{q_1} = a_q \text{ and } a_{q_1}^* = -a_q^* + \theta^*(2\pi - a_h) \tag{36}$$

**Result 4.2.** If $\varphi_q$ and $\varphi_{q_1}$ are the right closed Mannheim trajectory ruled surfaces, i.e., $\theta = \pi/2$, then the relationships between the real integral invariants of $\varphi_q$ and $\varphi_{q_1}$-CTRS are given as follows

$$\lambda_{q_1} = \lambda_h, \quad \ell_{q_1} = \ell_h + \theta^*\lambda_q \tag{37}$$

Then, the measure of spherical surface areas bounded by the spherical images of $\varphi_{q_1}$ and $\varphi_h$-CTRS are the same, i.e.,

$$a_{q_1} = a_h \text{ and } a_{q_1}^* = -\left(a_h^* + \theta^*(2\pi - a_q)\right) \tag{38}$$

**Result 4.3.** If $\theta^* = 0$, i.e., the generators $\vec{q}$ and $\vec{q}_1$ of the Mannheim offset surfaces intersect, then we have

$$\begin{cases} \lambda_{q_1} = \lambda_q \cos\theta + \lambda_h \sin\theta \\ \ell_{q_1} = \ell_q \cos\theta + \ell_h \sin\theta \end{cases} \tag{39}$$

In this case, $\varphi_q$ and $\varphi_{q_1}$-CTRS are intersect along their striction lines. It means, their striction lines are the same.

Let now consider that what the condition for the developable Mannheim offset of a CTRS is. Let $\varphi_q$ and $\varphi_{q_1}$-CTRS be the Mannheim offset surfaces and let $\vec{\alpha}(s)$ and $\vec{\beta}(s_1)$ be the striction lines of $\varphi_q$ and $\varphi_{q_1}$-CTRS, respectively. Then, we can write

$$\vec{\beta}(s) = \vec{\alpha}(s) + \theta^* \vec{a}(s) \tag{40}$$

where $s$ is the arc-length of $\vec{\alpha}(s)$. Assume that $\varphi_q$-CTRS is developable. Then from (3) and (6) we have

$$\delta_q = \frac{\langle \cos\sigma\vec{q} + \sin\sigma\vec{a},\ \vec{q} \times k_1\vec{h} \rangle}{\langle k_1\vec{h}, k_1\vec{h} \rangle} = \frac{\sin\sigma}{k_1} = 0 \tag{41}$$

Then we have $\sigma = 0$. Thus, from (6)

$$\frac{d\vec{\alpha}}{ds} = \vec{q} \tag{42}$$

Hence, along the striction line $\vec{\alpha}(s)$, the orthogonal frame $\{\vec{q}, \vec{h}, \vec{a}\}$ coincides with the Frenet frame $\{\vec{T}, \vec{N}, \vec{B}\}$ and the differential forms $k_1$ and $k_2$ turn into the curvature $\kappa_\alpha$ and torsion $\tau_\alpha$ of the striction line $\alpha(s)$, respectively. Then with the aid of (5), (40) and (42) we have

$$\frac{d\vec{\beta}}{ds} = \vec{q} - \theta^* \tau_\alpha \vec{h}. \tag{43}$$

On the other hand from (5) and (28) we obtain

$$\frac{d\vec{q}_1}{ds} = -\left(\frac{d\theta}{ds} + \kappa_\alpha\right)\sin\theta\vec{q} + \left(\frac{d\theta}{ds} + \kappa_\alpha\right)\cos\theta\vec{h} + \tau_\alpha \sin\theta\vec{a} \tag{44}$$

By using (29) and the fact that $k_1 = \kappa_\alpha$, from (44) we have

$$\frac{d\vec{q}_1}{ds} = \tau_\alpha \sin\theta\vec{a} \tag{45}$$

From (43) and (45) we have

$$\delta_{q_1} = \frac{\langle d\vec{\beta},\ \vec{q}_1 \times d\vec{q}_1 \rangle}{\langle d\vec{q}_1, d\vec{q}_1 \rangle} = \frac{\sin\theta + \theta^* \tau_\alpha \cos\theta}{\tau_\alpha \sin\theta}. \tag{46}$$

Thus, from (42) and (46) it can be stated that if the Mannheim offsets of $\varphi_q$ and $\varphi_{q_1}$ ruled surfaces are developable then the following relationship holds

$$\sin\theta + \theta^* \tau_\alpha \cos\theta = 0 \tag{47}$$

If (47) holds, along the striction line $\beta(s_1)$, the orthogonal frame $\{\vec{q}_1, \vec{h}_1, \vec{a}_1\}$ coincides with the Frenet frame $\{\vec{T}_1, \vec{N}_1, \vec{B}_1\}$. Thus, the following theorem may be given.

***Theorem 4.3.*** *If $\varphi_q$ and $\varphi_{q_1}$ are Mannheim offsets of developable trajectory ruled surfaces then their striction lines are Mannheim partner curves.*

From (47) we can give the following special cases

**Result 4.4.** $\theta = 0$, i.e., the Mannheim offsets of $\varphi_q$ and $\varphi_{q_1}$ developable trajectory ruled surfaces are oriented
  ⇔ Their generator are coincident, i.e., $\theta^* = 0$.
  ⇔ The Mannheim offsets of $\varphi_q$ and $\varphi_{q_1}$ developable trajectory ruled surfaces are coincident.

**Result 4.5.** $\theta = \pi/4$ ⇔ there is a relationship between the torsion of $\alpha(s)$ and offset distance as follows

$$\tau_\alpha = -\frac{1}{\theta^*} \tag{48}$$

If $\varphi_q$-CTRS is developable then from the equations (8), (28) and (43) the pitch $\ell_{q_1}$ of $\varphi_{q_1}$-CTRS is

$$\ell_{q_1} = -\oint (\cos\theta - \theta^* \tau_\alpha \cos\theta)ds$$

Then we can give the following result:

**Result 4.6.** If $\varphi_q$-CTRS is developable then the relation between the pitch $\ell_{q_1}$ of $\varphi_{q_1}$-CTRS and the torsion of sitriction line $\alpha(s)$ of $\varphi_q$-CTRS is given by

$$\ell_{q_1} = -\oint (\cos\theta - \theta^* \tau_\alpha \cos\theta)ds \tag{49}$$

Let now consider the area of projections of Mannheim offsets. With the aid of [5] the dual area vectors of the spherical images of $\varphi_q$ and $\varphi_{q_1}$ Mannheim offsets are

$$\begin{cases} \tilde{w}_q = \tilde{d} + \overline{\wedge}_q \tilde{q} \\ \tilde{w}_{q_1} = \tilde{d} + \overline{\wedge}_{q_1} \tilde{q}_1 \end{cases} \tag{50}$$

respectively. Then, the dual area of projection of the spherical image of $\varphi_{q_1}$-CTRS in direction $q$, generators of $\varphi_q$-offsets, is

$$2\overline{f}_{\tilde{q}_1,\tilde{q}} = \langle \tilde{w}_{q_1}, \tilde{q} \rangle = -\overline{\wedge}_q + \overline{\wedge}_{q_1} \cos\overline{\theta} \tag{51}$$

Separating (51) into real and dual parts we have the following theorem

***Theorem 4.4.*** *The relationships between the area of projections of spherical images of the Mannheim offsets and their integral invariants are given as follows*

$$\begin{cases} 2f_{q_1,q} = -\lambda_q + \lambda_{q_1} \cos\theta, \\ 2f^*_{q_1,q} = \ell_q - \ell_{q_1} \cos\theta - \lambda_{q_1}\theta^* \sin\theta, \end{cases} \tag{52}$$

**Result 4.7.** If $\varphi_q$ and $\varphi_{q_1}$-CTRS are the oriented surfaces, i.e., $\theta = 0$, then from (52) we have

$$2f_{q_1,q} = -\lambda_q + \lambda_{q_1}, \quad 2f^*_{q_1,q} = \ell_q - \ell_{q_1} \tag{53}$$

**Result 4.8.** If $\varphi_q$ and $\varphi_{q_1}$-CTRS are the right Mannheim offsets, i.e., $\theta = \pi/2$, then from (52) we have

$$2f_{q_1,q} = -\lambda_q, \quad 2f^*_{q_1,q} = \ell_q - \lambda_{q_1}\theta^*, \tag{54}$$

Similarly, the dual area of projection of spherical image of $\varphi_{q_1}$-CTRS in the direction $\vec{h}$ is

$$2\overline{f}_{q_1,h} = \langle \tilde{w}_{q_1}, \tilde{h} \rangle = -\overline{\wedge}_h + \overline{\wedge}_{q_1} \sin\overline{\theta} \tag{55}$$

Separating (51) into real and dual parts we have the following

**Result 4.9.** If $\varphi_q$ and $\varphi_{q_1}$-CTRS are Mannheim offsets then we have

$$\begin{cases} 2f_{q_1,h} = -\lambda_h + \lambda_{q_1} \sin\theta, \\ 2f^*_{q_1,h} = \ell_h - \ell_{q_1} \sin\theta + \lambda_{q_1}\theta^* \cos\theta, \end{cases} \tag{56}$$

**Result 4.10.** If $\varphi_q$ and $\varphi_{q_1}$-CTRS are the oriented Mannheim surfaces, i.e., $\theta = 0$, then from (56) we have

$$2f_{q_1,h} = -\lambda_h, \quad 2f^*_{q_1,h} = \ell_h + \lambda_{q_1}\theta^* \tag{57}$$

**Result 4.11.** If $\varphi_q$ and $\varphi_{q_1}$-CTRS are right Mannheim offsets, i.e., $\theta = \pi/2$, then from (56) we have

$$2f_{q_1,h} = -\lambda_h + \lambda_{q_1}, \quad 2f^*_{q_1,h} = \ell_h - \ell_{q_1} \tag{58}$$

Similarly, the dual area of projection of spherical image of $\varphi_{q_1}$-CTRS in the direction $\vec{a}$ is

$$2\bar{f}_{q_1,a} = \langle \tilde{w}_{q_1}, \tilde{a} \rangle = -\overline{\wedge}_a = -\overline{\wedge}_{h_1} = 0 \tag{59}$$

Separating (59) into real and dual parts we have the following

**Result 4.12.** If $\varphi_q$ and $\varphi_{q_1}$-CTRS are Mannheim offsets then we have

$$f_{q_1,a} = -\lambda_{h_1} = -\lambda_a = 0, \quad f^*_{q_1,a} = \ell_{h_1} = \ell_a = 0 \tag{60}$$

**Example 4.1.** Let consider the cone surface given by
$$\varphi(u,v) = (0,1,0) + v(\cos u, 1+\sin u, 1) \tag{61}$$
(Fig 1), which corresponds to dual curve

$$\tilde{q}(u) = \frac{1}{\sqrt{2}}(\cos u, \sin u, 1) + \varepsilon \frac{1}{\sqrt{2}}(1, 0, -\cos u).$$

Some Mannheim offsets of the surface in (61) can be given as follows:

**i)** The oriented Mannheim offset of (61) with dual offset angle $\bar{\theta} = 0 + \varepsilon\sqrt{2}$ is given by

$$\varphi_1(u,v) = (-\cos u, 1-\sin u, 1) + v\frac{1}{\sqrt{2}}(\cos u, \sin u, 1) \tag{62}$$

which is also a cone (Fig. 4).

**ii)** The right Mannheim offset of (61) with dual offset angle $\bar{\theta} = \pi/2 + \varepsilon\sqrt{2}u$ is given by
$$\varphi_2(u,v) = (-u\cos u, 1-u\sin u, u) + v(-\sin u, \cos u, 0) \tag{63}$$
which represents a helicoids (Fig. 5).

**iii)** The Mannheim offset of (61) with dual offset angle $\bar{\theta} = \pi/3 + \varepsilon\sqrt{2}$ is

$$\varphi_3(u,v) = (-\cos u, 1-\sin u, 1) + v(\frac{1}{2\sqrt{2}}\cos u - \frac{\sqrt{3}}{2}\sin u, \frac{1}{2\sqrt{2}}\sin u + \frac{\sqrt{3}}{2}\sin u, \frac{1}{2\sqrt{2}}) \tag{64}$$

which is a hyperboloid of one sheet (Fig.6).

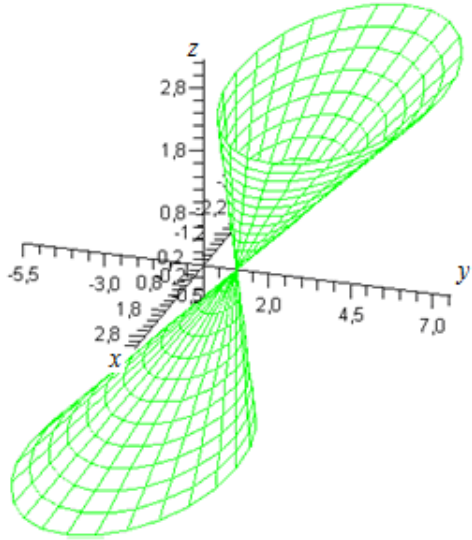

**Fig. 3.** The surface $\varphi(u,v)$.

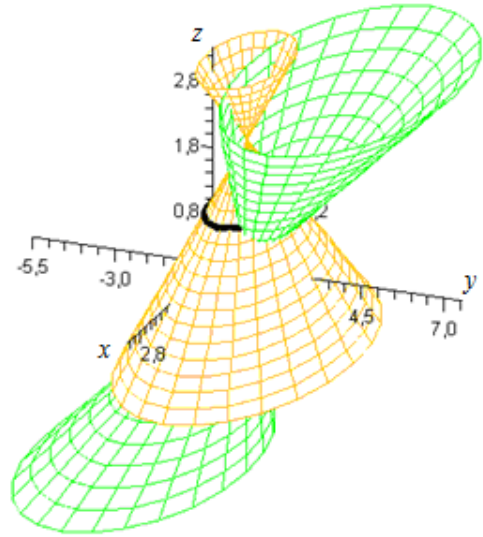

**Fig. 4.** Mannheim offset $\varphi_1(u,v)$ of $\varphi(u,v)$ with dual offset angle $\bar{\theta} = 0 + \varepsilon\sqrt{2}$.

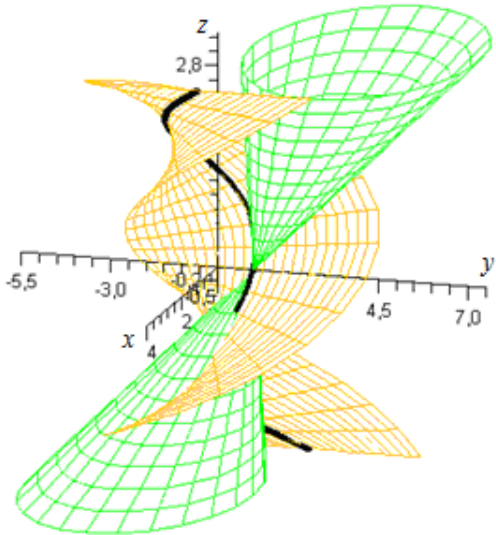

**Fig. 5.** Mannheim offset $\varphi_2(u,v)$ of $\varphi(u,v)$ with dual offset angle $\bar{\theta} = \pi/2 + \varepsilon\sqrt{2}u$.

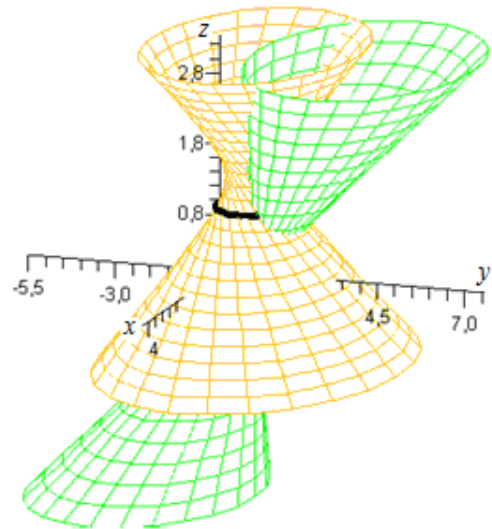

**Fig. 6.** Mannheim offset $\varphi_3(u,v)$ of $\varphi(u,v)$ with dual offset angle $\bar{\theta} = \pi/3 + \varepsilon\sqrt{2}$.

In Figures 4, 5 and 6, the curves rendered in black are the striction lines of the Mannheim offsets of the surface $\varphi_i(u,v),\ (1 \leq i \leq 3)$.

## 4. Conclusion

In this paper, we give the characterizations of Mannheim offsets of ruled surfaces. We find new relations between the invariants of Mannheim offsets of ruled surfaces. Furthermore, we show that the striction lines of the Mannheim offsets of the developable ruled surface are Mannheim partner curves. By using the similar methods given for Bertrand offsets, Mannheim offsets of the ruled surfaces can be used in CAGD.